\documentstyle[12pt]{article}
\begin{document}
\begin{titlepage}
\hspace*{12.0cm}HU-TFT-93-28
\vskip 2.0 cm

\begin{center}
{\Large\bf q-Supersymmetric Generalization of von Neumann's Theorem}
\vskip 1.0 cm
by
\vskip 0.5 cm
{\bf M. Chaichian}$^1$,
{\bf R. Gonzalez
Felipe}$^{1,}$\renewcommand{\thefootnote}{*)}\footnote{ICSC-World
Laboratory; On leave of absence from Grupo de F\'isica
Te\'orica, Instituto de Cibern\'etica, Matem\'atica y Fisica, Academia de
Ciencias de Cuba, Calle E No. 309, Vedado, La Habana 4, Cuba.} {\bf and
P. Pre\v{s}najder}$^{2,}$\renewcommand{\thefootnote}{**)}\footnote{Permanent
address:
Department of Theoretical Physics, Comenius University, Milynsk\'a dolina
F2, CS-84215 Bratislava, Slovakia.}
\vskip 1.0 cm
$^1$High Energy Physics Laboratory, Department of Physics,
P.O. Box 9 (Siltavuorenpenger 20 C), SF-00014 University of Helsinki, Finland
\vskip 1.0 cm
$^2$Research Institute for Theoretical Physics, University of Helsinki,
P.O. Box 9 (Siltavuorenpenger 20 C), SF-00014 University of Helsinki, Finland
\vskip 2.0 cm
{\bf May 1993}

\end{center}

\vskip 3.5 cm

{\bf Abstract}

Assuming that there exist operators which form an irreducible representation
of the $q$-superoscillator algebra, it is proved that any two such
representations are equivalent, related by a uniquely determined superunitary
transformation. This provides with a $q$-supersymmetric generalization of the
well-known uniqueness theorem of von Neumann for any finite number of degrees
of freedom.
\end{titlepage}

\section{Introduction}

\hspace*{1.0cm} In the last few years quantum deformations of Lie groups and
Lie
algebras have
found several applications in mathematics and theoretical physics (see e.g.
refs. \cite{fad}-\cite{proc1}). These deformations have been subsequently
extended to
supergroups and superalgebras \cite{yu}-\cite{chaic}. In particular, the
bosonic and fermionic
$q$-oscillators \cite{mac}, \cite{bie}, \cite{chaic} have been used for the
realization of different
quantum Lie algebras \cite{mac}-\cite{hay} and quantum superalgebras
\cite{chaic},\cite{chaku},\cite{chakulu}.

A natural question then arises, concerning the relation between different
irreducible representations of the $q$-deformed algebras. It is well known
that in the case of the classical bosonic and fermionic Heisenberg algebras
for harmonic oscillators, this problem is solved by the von Neumann's theorem
(see, e.g. \cite{ber},\cite{put}), which states that irreducible
representations of the
bosonic (fermionic) algebra are unitarily equivalent to each other. Similar
results also hold for irreducible operator representations of Lie superalgebras
\cite{gros}. Recently, it was proved that the analogue of von Neumann's
theorem is also valid in the case of $q$-oscillator algebras \cite{chamnat}.

In this letter we shall extend the results of \cite{gros} and formulate a
quantum
supersymmetric generalization of von Neumann's theorem for irreducible
representations of $q$-deformed superalgebras. We start from a supercovariant
system of $q$-oscillators \cite{chakulu}, which are covariant under the
coaction of
a supergroup, $SU_q(n|m)$. The latters present the extension of the covariant
system of $q$-oscillators proposed in \cite{pusz},\cite{wpus}. Assuming
suitable domain
properties as in \cite{gros}, we prove that any two representations of the
$q$-deformed
superalgebra are connected by a unique {\it superunitary} transformation.
A similar result is also proved to be valid for any finite number, $n$, of
bosonic and fermionic independent $q$-oscillators. We present the explicit
form of the superunitary transformation operator for the cases $n=1$ and 2.
\vskip 0.5 cm

\section{von Neumann's theorem and its extension to $q$-deformed algebras}

\hspace*{1.0cm} We start by recalling the classical von Neumann's theorem
\cite{ber},\cite{put}. Let
$b,b^+$ and $b',b'^+$ be two irreducible representations of the Heisenberg
algebra,
$$bb^+-b^+b=1\ ,
\eqno{(1)}$$
in the Hilbert spaces $H$ and $H'$, respectively, and assume that there exist
vectors $|0>$ in $H$ and $|0'>$ in $H'$, such that $b|0>=0,\ b'|0'>=0$. Then
there exists a unitary operator $U$ such that
$$
\begin{array}{l}
b'=UbU^+\ ,\ b'^+=Ub'U^+\ ,\\
\ \\
UU^+=U^+U=1\ .
\end{array}\eqno{(2)}$$
A similar theorem also holds in the case of the fermionic algebra
\cite{ber},\cite{put},
$$cc^++c^+c=1\ .
\eqno{(3)}$$
\indent The above theorem can be formulated in a slightly different form,
which we shall utilize further:

Let $b,b^+$ and $b',b'^+$ satisfy all the conditions of von Neumann's theorem
and let us define the operators $a,a^+$ and $a',a'^+$ as
$$
\begin{array}{lll}
a=\varphi(N)b&,&a^+=b^+\varphi^+(N)\ ,\\
a'=\varphi(N')b'&,&a'^+=b'^+\varphi^+(N')\ ,\end{array}
\eqno{(4)}$$
where
$$N=b^+b\ ,\ N'=b'^+b'\ ,
\eqno{(5)}$$
are the number operators and $\varphi$ is a well-behaved function.

Then for the representations $a,a^+$ and $a',a'^+$ in some Hilbert (sub)spaces
$H_a$ and $H'_a\ (H_a\subseteq H,H'_a\subseteq H')$, respectively,
von Neumann's theorem also holds, i.e., $a,a^+$ and $a',a'^+$ are
irreducible representations, there exist vectors $|0>$ in $H_a$ and $|0'>$ in
$H'_a$ such that $a|0>=0,\ a'|0'>=0$ and there exists a unitary operator $U$
such
that
$$ a'=UaU^+\ ,\ a'^+=Ua^+U^+\ ,\ N'=UNU^+\ .
\eqno{(6)}$$
\indent To prove (6) we notice that since $b,b^+$ and $b',b'^+$ satisfy the
conditions of von
Neumann's theorem, it follows that there exists a unitary operator $U$ such
that relations (2) hold. From (2) and (5) it follows that
$$N'=b'^+b'=Ub^+U^+UbU^+=UNU^+\ ,$$
and the same relation is also valid for the function $\varphi (N)$,
$$
\varphi(N')=U\varphi(N)U^+\ .
\eqno{(7)}$$
Now, from the definition of the operators $a,a^+$ (eqs. (4)) we have
$$UaU^+=U\varphi(N)bU^+=U\varphi(N)U^+UbU^+=\varphi(N')b'=a'\ ,$$
and, similarly,
$$Ua^+U^+=a'^+\ .$$
Relations (6) are thus proved.

It is also clear that in the spaces $H_a$ and $H'_a$ there exist vectors
$|0>$ and $|0'>$ such that $a|0>=0,\ a'|0'>=0$.
Finally, $a,a^+$ and $a',a'^+$ are two irreducible representations in the
Hilbert spaces $H_a,\ H'_a$, respectively. This follows from the definitions
(4)
and the fact that $b,b^+$ and $b',b'^+$ are irreducible representations in
$H,\ H'$.

The same conclusion also holds obviously for the case of fermionic
oscillators.

{}From the above theorem, it follows, in particular, that any two irreducible
representations of the $q$-deformed bosonic oscillator algebra
\cite{mac},\cite{bie},\cite{chaic}
$$
\begin{array}{l}aa^+-qa^+a=q^{-N}\ ,\\
$$[N,a]=-a\ ,\ [N,a^+]=a^+\ ,\end{array}
\eqno{(8)}$$
are connected by a unique unitary transformation. Indeed, by taking
$\varphi(N)=\sqrt{\frac{[N+1]}{N+1}}$ in (4), where $[n]=\frac{q^n-q^{-n}}
{q-q^{-1}}$, we obtain the $q$-oscillator algebra (8), provided that
$a^+a=[N],\ aa^+=[N+1]$, i.e. in the Fock space of (8) \cite{kuda}.
For the exceptional values of $q$ being the $m$-th root of unity,
$q=e^{\pm i\pi/m},[m]=0$ and thus $\varphi(m-1)=0$. In such a case the
Hilbert space $H_a$ becomes finite, $m$-dimensional.
A similar statement is also valid in the case of the $q$-deformed fermionic
algebra \cite{chaic}
$$
\begin{array}{l}
ff^++qf^+f=q^M\ ,\\
$$[M,f]=-f\ ,\ [M,f^+]=f^+\ ,\end{array}\eqno{(9)}$$
which can be obtained from the algebra (3) by means of the change of operators
$$f=q^{\frac{M}{2}}c\ ,\ f^+=c^+q^{\frac{M}{2}}\ .$$

Let us remark that the above statements were proved in ref. \cite{chamnat} by
following
a proof in the same line as the original von Neumann's theorem for usual
(nondeformed) oscillators and valid also for generic values of $q$ including
the
exceptional values of $m$-th root of unity. The same theorem
is also valid in the case of the other $q$-fermionic algebra \cite{parth} given
by the
commutation relations
$$ff^++qf^+f=q^{-M}\ ,$$
$$[M,f]=-f\ ,\ [M,f^+]=f^+\ .$$
In fact, we can obtain the latters from the bosonic algebra (1) by means of the
transformation,
$$
f=\sqrt{\frac{[M+1]^f}{M+1}}b\ ,\ f^+=b^+\sqrt{\frac{[M+1]^f}{M+1}}\ ,$$
where $[n]^f=\frac{q^{-n}-(-1)^nq^n}{q+q^{-1}}$, since in the Fock space we
have the relations $f^+f=[M]^f\ ,\ ff^+=[M+1]^f$\ .
\vskip 1.0 cm

\section{$q$-Supersymmetric von Neumann's Theorem}

\hspace*{1.0cm}In this section we formulate our main result, namely, we prove
that an irreducible representation of the $q$-deformed superalgebras is, up to
a
superunitary transformation, unique, in the following sense:

{\bf Theorem (i)}

Let $Z=\{B,B^+,F,F^+\}$ be an irreducible operator family, which satisfies the
$q$-deformed superalgebra [11]:
$$
\begin{array}{lll}
BF&=&qFB\ ,\ B^+F^+=q^{-1}F^+B^+\ ,\end{array}
\eqno{(10{\rm a})}$$
$$
\begin{array}{lll}
BF^+&=&q^{-1}F^+B\ ,\ B^+F=qFB^+\ ,\end{array}
\eqno{(10{\rm b})}$$
$$
\begin{array}{lll}
F^2&=&(F^+)^2=0\ ,\hspace{2.5cm}\end{array}
\eqno{(10{\rm c})}$$
$$
\begin{array}{lll}
BB^+-q^{-2}B^+B&=&1+(q^{-2}-1)F^+F\ ,\end{array}
\eqno{(10{\rm d})}$$
$$
\begin{array}{llll}FF^++F^+F=1;&\end{array}
\eqno{(10{\rm e})}$$
where $q$ is a real number, $F,F^+$ are bounded operators on the separable
Hilbert space $H$ and $B,B^+$ are densely defined closed operators in $H$.
Let $\theta$ be a Grassmann variable, such that
$$\{\theta, F\}=\{\theta, F^+\}=\theta^2=0\ ,
\eqno{(11{\rm a})}$$
$$[\theta,B]=[\theta,B^+]=0\ .\hspace{1.0cm}
\eqno{(11{\rm b})}$$
Let $D$ and $G$ be densely defined closed linear operators and define on a
suitable domain \cite{gros}
$$B'=B+\theta D\ ,\ \ B'^+=B^+-\theta D^+\ ,
\eqno{(12)}$$
$$F'=F+\theta G\ ,\ \ F'^+=F^++\theta G^+\ ,$$
\
$$G=G_{00}(B,B^+)+G_{11}(B,B^+)F^+F\ ,$$
$$D=D_{10}(B,B^+)F^++D_{01}(B,B^+)F\ ,\eqno{(13)}$$

\noindent where $G$ and $D$ are assumed to be even and odd Grassmann elements,
respectively. Assume that the  operator family $Z'=\{B',B'^+,F',F'^+\}$ also
fulfills the algebra (10) on a suitable domain of definition.

Then, under the above conditions, there exists a uniquely determined
self-adjoint odd operator $A$, such that
$$G=\{A,F\}\ ,\ D=[A,B]\ ,
\eqno{(14)}$$
$$A=F^+G_{00}(B,B^+)\ +\ G^+_{00}(B,B^+)F=A^+\ ,
\eqno{(15)}$$
and the transformation (12) is implemented by the superunitary operator
$e^{\theta A}$ such that
$$e^{\theta A}Be^{-\theta A}=B+\theta [A,B]=B'\ ,$$
$$e^{\theta A}Fe^{-\theta A}=F+\theta\{A,F\}=F'\ .
\eqno{(16)}$$
Under the conditions of the above theorem, we have
$$G_{11}(B,B^+)=G_{00}(qB,qB^+)-G_{00}(B,B^+)\ ,
\eqno{(17{\rm a})}$$
$$
\begin{array}{lll}
D_{10}(B,B^+)&=&qG_{00}(qB,qB^+)B-BG_{00}(qB,qB^+)\ ,\\
\ \\
D_{01}(B,B^+)&=&q^{-1}G^+_{00}(B,B^+)B-BG^+_{00}(B,B^+)\ .\end{array}
\eqno{(17{\rm b})}$$
The transformation $e^{\theta A}=1+\theta A$ is called
{\it superunitary} if $A$ is odd and self-adjoint. In particular, from this
definition it follows that $\{\theta, A\}=0$ and the latter is fulfilled
only when $\theta$ satisfies the commutation relations (11).
\pagebreak

{\bf Remark:}

The Grassmann variable $\theta$, being on the $q$-plane of the
oscillators given in (10), can in general have the $q$-commutation relations
with the elements of the algebra as:
$$
\begin{array}{ll}
\theta_qF+pF\theta_q=0\ ,& \theta_qF^++p^{-1}F^+\theta_q=0\ ,\\
\theta_qB-rB\theta_q=0\ ,& \theta_qB^+-r^{-1}B^+\theta_q=0\ ,\end{array}
$$
with $p$ and $r$ real numbers
\renewcommand{\thefootnote}{*}\footnote{These $q$-commutation relations,
because of associativity, are of course compatible with the supersymmetric
$q$-Jacobi identities
$$[[A,B\}_{(q_3,q^{-1}_3)}, C \}_{(\frac{q_1}{q_2},\frac{q_2}{q_1})}
+(-1)^{\eta_A(\eta_B+\eta_C)}[[B,C \}_{(q_1,q^{-1}_1)},A \}_{(\frac{q_2}{q_3},
\frac{q_3}{q_2})}$$
$$+(-1)^{\eta_C(\eta_A+\eta_B)}[[C,A\}_{(q_2,q^{-1}_2)},B\}
_{(\frac{q_3}{q_1},\frac{q_1}{q_3})}=0\ ,$$
where $\eta_Z=1$ if $Z$ is odd, $\eta_Z=0$ if $Z$ is even and $[A,B\}_{(p,q)}
\equiv pAB\mp qBA$, where we take the plus sign when both $A$ and $B$ are odd,
otherwise we take the minus sign. The above expression represents the most
general form of the supersymmetric $q$-Jacobi identities, which includes three
arbitrary complex parameters $q_1,q_2$ and $q_3$.}. If however, we use the
above $q$-commutation relations with $p$ and $r$, instead of commutation
relations (11), from the requirement of oddness of operator $A$ and the
superunitarity condition $(e^{\theta_qA})^+=e^{-\theta_qA}$, i.e.
$A^+\theta_q=-\theta_qA$, one arrives at the following "$q$-self adjointness"
condition,
$$
A^+(B,B^+,F,F^+)=A(rB,r^{-1}B^+,pF,p^{-1}F^+)\ ,$$
instead of the usual one. The latter restriction seems rather unnatural, since
a usual self-adjoint operator $A$ on the original algebra of oscillators (10)
now acquires restrictions by the inclusion of an additional auxiliary
Grassmann element $\theta_q$. Thus for our purpose we can restrict ourselves
to the case $p=r=1$, i.e., to the commutation relations (11).

{\bf Proof of the theorem:}

The proof is actually quite similar to the one given
in \cite{gros} for the classical Lie superalgebras.

{}From the relations $F^2=F'^2=0$ and using (11a) we have, according to (12),
$F'^2=\theta [G,F]=0$, from which we obtain
$$[G,F]=0\ .
\eqno{(18)}$$
Relations (10a), (10b) imply that for any function $g(B,B^+)$ it holds
$$
\begin{array}{lll}
g(B,B^+)F&=&Fg(qB,qB^+)\ ,\\
g(B,B^+)F^+&=&F^+g(q^{-1}B,q^{-1}B^+)\ .\end{array}
\eqno{(19)}$$
With the ansatz (13) and using (10e) and (19) we find
$$[G,F]=\{G_{00}(B,B^+)-G_{00}(q^{-1}B,q^{-1}B^+)-G_{11}(q^{-1}B,q^{-1}B^+)\}F\
{}.
\eqno{(20)}$$
Comparing (18) and (20) we get
$$G_{11}(B,B^+)=G_{00}(qB,qB^+)-G_{00}(B,B^+)\ ,
\eqno{(21)}$$
and substituting it into (13) we obtain
$$G=G_{00}(B,B^+)FF^++G_{00}(qB,qB^+)F^+F\ .
\eqno{(22)}$$
Let us now take $G$, in the required form (14), as
$$G=\{A,F\}\ ,
\eqno{(23)}$$
where $A$ is odd. Assuming that $A$ is self-adjoint, i.e. $A^+=A$, we can
write
$$A=F^+\alpha (B,B^+) + \alpha^+(B,B^+)F
\eqno{(24)}$$
and we have from (23),
$$G=\alpha(qB,qB^+)F^+F+\alpha(B,B^+)FF^+
\eqno{(25)}$$
Comparing (22) and (25), we find
$$\alpha(B,B^+)=G_{00}(B,B^+)\ ,$$
and by substituting it back into (24), the relation (15) is obtained.

Now to find $D$, we use relations (10a) and (11). We have
$$B'F'-qF'B'=\theta\{ DF+BG-qGB+qFD\}=0\ .
\eqno{(26)}$$
Substituting (23) into (26) and using (10a), we will obtain then that
$$\{D-[A,B]\}F+qF\{D-[A,B]\}=0\ .
\eqno{(27{\rm a})}$$
Similarly, from (10b), (11) and $G^+=\{A,F^+\}$ we obtain
$$\{D-[A,B]\}F^++q^{-1}F^+\{D-[A,B]\}=0\ .
\eqno{(27{\rm b})}$$
Thus, from eqs. (27) we conclude that
$$D=[A,B]\ .
\eqno{(28)}$$
Finally, the operator $e^{\theta A}$ is superunitary since $A$ is odd and
self-adjoint (see eq. (24)). Moreover, we have
$$e^{\theta A}Be^{-\theta A}=(1+\theta A)B(1-\theta A)=B+\theta[A,B]=B'\ ,$$
$$e^{\theta A}Fe^{-\theta A}=(1+\theta A)F(1-\theta A)=F+\theta\{A,F\}=F'\ ,$$
i.e. relations (16) hold. The proof of the relations
(17) is straighforward: eq. (17a) was already proved (see eq. (21)) and
relations (17b) follow from (13), (14), (15) and (19). This completes the
proof.
When $q=1$, we reproduce the results of \cite{gros}.

Let us remark that a pair of one bosonic, $b,b^+$, and one fermionic,
$f,f^+$, independent
$q$-oscillators (which satisfy the relations (30)), can also be introduced
by means of the transformation
\cite{chakulu}
$$
\begin{array}{lll}
\ b=q^{\frac{N}{2}+M}B&,&f=q^{\frac{M}{2}}F\ ,\\
\ [N,B]=-B&,&[N,B^+]=B^+\ ,\ N^+=N\ ,\\
\ [M,F]=-F&,&[M,F^+]=F^+\ ,\ M^+=M\ ,\end{array}
\eqno{(29)}$$
where $B$ and $F$ are the elements of supercovariant algebra (10).
In this case, the uniqueness
(up to a superunitary transformation) of any irreducible representation of
$b,b^+,N,f,f^+,M$ is given by the following theorem.

{\bf Theorem (ii)}

Let $Z=\{b,b^+,N,f,f^+,M\}$ be an irreducible operator family, which satisfies
the $q$-oscillator algebra
$$
\begin{array}{lll}
[b,f]=0\ ,\ [b,f^+]&=&0\ ,\end{array}
\eqno{(30{\rm a})}$$
$$
\begin{array}{lll}
f^2=(f^+)^2&=&0\ ,\hspace{1.0cm}\end{array}
\eqno{(30{\rm b})}$$
$$
\begin{array}{lll}
bb^+-q^{-1}b^+b&=&q^N\ ,\hspace{0.5cm}\end{array}
\eqno{(30{\rm c})}$$
$$
\begin{array}{lll}
ff^++qf^+f&=&q^M\ ,\hspace{1.0cm}\end{array}
\eqno{(30{\rm d})}$$
$$
\begin{array}{lll}
[N,b]=-b,\;[N,b^+]&=&b^+\ ,\end{array}
\eqno{(30{\rm e})}$$
$$
\begin{array}{lll}
[M,f]=-f,\;[M,f^+]&=&f^+\ ,\end{array}
\eqno{(30{\rm f})}$$
on a suitable domain of definition; $q$ is real. Let $\theta$ be a Grassmann
variable,
$$
\begin{array}{l}
\ \{\theta,f\}=\{\theta,f^+\}=\theta^2=0\ ,\\
\ \\
\ [\theta,b]=[\theta,b^+]=0\ .\end{array}
\eqno{(31)}$$
Define
$$
b'=b+\theta D,\ \ \ f'=f+\theta G,
$$
$$
D=D_{10}(b,b^+,M)f^++D_{01}(b,b^+,M)f\ ,
\eqno{(32)}
$$
$$
G=G_{00}(b,b^+,M)+G_{11}(b,b^+,M)f^+f\ ,$$
where $D$ and $G$ are the even and odd Grassmann elements, respectively,
and have the most general form as in (32).
Assume that the  operator family $Z'=\{b',b'^+,N',f',f'^+,M'\}$ also fulfills
the algebra (30) on a suitable domain of definition.

Then there exists a uniquely determined self-adjoint odd operator $A$, such
that
one can write
$$
G=\{A,f\}\ ,\ \ D=[A,b]\ ,
\eqno{(33)}$$
$$A=f^+\alpha(b,b^+,M)+\alpha^+(b,b^+,M)f=A^+
\eqno{(34)}$$
and the transformation (32) is implemented by the superunitary operator
$e^{\theta A}$ such that
$$e^{\theta A}be^{-\theta A}=b'\ ,$$
$$e^{\theta A}fe^{-\theta A}=f'\ .
\eqno{(35)}$$
Under the conditions of the above Theorem (ii), we have
$$G_{00}(b,b^+,M)=\alpha(b,b^+,M)q^M,\ G_{11}(b,b^+,M)=\alpha(b,b^+,M-1)-q
\alpha(b,b^+,M)
\eqno{(36{\rm a})}$$
and
$$D_{10}(b,b^+,M)=[\alpha(b,b^+,M-1),b],\ D_{01}=[\alpha^+(b,b^+,M),b].
\eqno{(36{\rm b})}$$
The proof is similar to the one of Theorem (i) and will not be given.

Let us notice that from the Theorem (ii), it also follows that
$$e^{\theta A}Ne^{-\theta A}=N+\theta[A,N]=N'\ ,$$
$$e^{\theta A}Me^{-\theta A}=M+\theta[A,M]=M'\ .
\eqno{(37)}$$
Indeed, relations (35) imply that for any functions $\varphi(bb^+,b^+b)$
and $\psi(ff^+,f^+f)$, one has
$$e^{\theta A}\varphi(bb^+,b^+b)e^{-\theta A}=\varphi(b'b'^+,b'^+b')\ ,$$
$$e^{\theta A}\psi(ff^+,f^+f)e^{-\theta A}=\psi(f'f'^+,f'^+f')\ ;$$
$\varphi=N$ and $\psi=M$ are just particular cases of these functions (see
eqs. (30c) and (30d) respectively).

Here we would like to mention that there exists a relation between the
transformations of
Theorem (i) and of Theorem (ii). Indeed, if one takes in (34)
$$\alpha(b,b^+,M)=q^{-\frac{M}{2}}G_{00}(B,B^+)\ ,
\eqno{(38)}$$
where $B=q^{-\frac{N}{2}-M}b,\ F=q^{-\frac{M}{2}}f$ according to (29), it is
straighforward to show that the superunitary transformation generated by (34),
(35) with $\alpha$ given in (38), corresponds to the superunitary
transformation (15), (16) of Theorem (i), and therefore eqs. (17) are
satisfied.

Also it will be interesting to find a direct relation between the
superunitary
transformations of the usual (undeformed) and of the $q$-deformed cases
obtained here, in the same way as such a relation exists for the
nonsupersymmetric $q$-oscillators treated in section 2.

To summarize, we have shown that the irreducible representations of the
$q$-superoscillator algebra are equivalent and are related by a unique
superunitary transformation. The $q$-supersymmetric generalization of von
Neumann's theorem, presented above, is ($N=1$ supersymmetry) for only one
bosonic and one fermionic degrees of freedom. Our results can be extended
to the supersymmetric case with any number of bosonic and fermionic degrees
of freedom. The theorem now can be formulated as follows:

{\bf Theorem (finite degrees of freedom):}

Let $Z=\{b_k,b^+_k,N_k,f_k,f^+_k,M_k;\ k=1,...,n\}$ be an irreducible operator
set, which satisfies the $q$-oscillator algebra
$$
\begin{array}{ll}
b_kb^+_k-q^{-1}b^+_kb_k=q^{N_k}\ ,&\\
\ \\
f_kf^+_k+q\;f^+_kf_k=q^{M_k}\ ,&\end{array}$$
$$[N_k,b_k]=-b_k\ ,\ [N_k,b^+_k]=b^+_k\ ,\eqno{(39)}$$
$$[M_k,f_k]=-f_k\ ,\ [M_k,f^+_k]=f^+_k\ ,$$
with all the other (anti)commutation relations vanishing (i.e., independent
system of $q$-oscillators); $q$ is real.

Assume now that another set of operators $Z'=\{b'_k,b'^+_k,N'_k,f'_k,f'^+_k,
M'_k;\ k=1,...n\}$ also satisfies the same algebra (39). Then the two sets
$Z$ and $Z'$ are equivalent, related by a unique superunitary transformation
such that
$$
\begin{array}{ll}
b'_k=Ub_kU^+\ ,& f'_k=Uf_kU^+\ ,\\
U=e^{\theta A}\ ,&\ \end{array}
\eqno{(40)}$$
with the same relation between the remaining elements of the two sets. In (40)
$\theta$ is a Grassmann variable satisfying the (anti)commutation relations
(31) for all the $b_k,b^+_k,f_k,f^+_k$ and $A$ is a self-adjoint odd operator.

The proof of the theorem can be most easily performed by the method of
induction in the number of degrees of freedom, $n$. In addition to the
explicit form of the superunitary operator for $n=1$ given in Theorem (ii), we
present below explicitly the formulae for the case $n=2$.

Define
$$b'_k=b_k+\theta D_k\ ,\ f'_k=f_k+\theta G_k\ ,\ \ k=1,2\ ,\eqno{(41)}$$
where $G_k$ have the most general form
$$G_k=G^0_k+G^1_kf^+_1f_1+G^2_kf^+_2f_2+G^{12}_kf^+_1f_1f^+_2f_2$$
$$+H_kf_1f_2+H'_kf^+_1f^+_2+K_kf_1f^+_2+K'_kf^+_1f_2\ ;\eqno{(42)}$$
all the coefficients in (42) can depend on $b_{\ell},b^+_\ell,N_\ell,M_\ell\ \
(\ell=1,2).$\\
{}From the anticommutation relations $\{f'_k,f'_\ell\}=\{f'^+_k,f'^+_\ell\}=
\{f'_1,f'^+_2\}=\{f'^+_1,f'_2\}=0$, we obtain that $H'_k=K'_k=0$ and that
the coefficients $G^k_k,G^{12}_k,H_k,K_k\ (k=1,2)$ can be expressed in terms
of $G^0_1,G^0_2,G^1_2,G^2_1$ and their hermitian conjugates. Now if we write
$G_k$ in the form $G_k=\{A,f_k\}$, then it is straightforward to show that
there exists a unique self-adjoint odd operator $A$ given by
$$
A=\sum\limits_kf^+_kG^0_kq^{-M_k}+\sum\limits_{k\neq\ell}f^+_kf^+_\ell f_\ell\,
G^\ell_k\,q^{-M_k}+\ h.c.\ ,$$
and from the commutativity between the $b'_k$ and $f'_\ell$, we obtain the
required form $D_k=[A,b_k]$.

\vskip 2.0 cm

{\bf Acknowledgements}

It is our pleasure to thank A. Demichev for useful discussions and several
clarifying remarks.
R.G.F. would like to thank ICSC-World Laboratory for financial support. P.P. is
grateful to the Research Institute for Theoretical Physics, University of
Helsinki, for the hospitality.

\end{document}